\documentclass{birkmult}

\newtheorem{thm}{Theorem}
 
 \newtheorem{lem}{Lemma}
 
 \theoremstyle{definition}
 \newtheorem{defn}{Definition}
 \theoremstyle{remark}

\usepackage{amsfonts,amsmath,amssymb,graphics}

\newcommand{\C}{\mathbb{C}}
\newcommand{\N}{\mathbb{N}}

\newcommand{\R}{\mathbb{R}}
\newcommand{\Z}{\mathbb{Z}}

\newcommand{\myfunc}{\phi}

\newcommand{\I}{\mathrm{i}}
\newcommand{\baromega}{\overline{\omega}}
\newcommand{\Wr}[2]{\mathcal{W}\left\{#1,#2\right\}}

\newcommand{\Ai}{\mathrm{Ai}}
\newcommand{\Bi}{\mathrm{Bi}}
\newcommand{\Ei}{\mathrm{Ei}}

\newcommand{\Uplus}{U^\oplus}
\newcommand{\Uminus}{U^\ominus}
\newcommand{\sequence}{S}

\newcommand{\notthis}[1]{}
\newcommand{\Real}{\mathrm{Re}\,}
\newcommand{\Imag}{\mathrm{Im}\,}
\newcommand{\tr}{\mathrm{tr}\,}

\newcommand{\diff}[1]{\left[#1\right]_-}

\newcommand{\half}{\frac{1}{2}}

\newcommand{\shalf}{{\textstyle\frac{1}{2}}}

\begin{document}

\author{R. J. Martin}

\address{
(Affiliated to) \\ Department of Mathematics\\ University College London \\ London WC1E 6BT, U.K. 
}

\email{rmartin@ahl.com}

\author{M. J. Kearney}
\address{ Faculty of Engineering and Physical Sciences \\ University of Surrey\\
Guildford, Surrey GU2 7XH, U.K. }

\email{m.j.kearney@surrey.ac.uk}


\title{\bf An exactly solvable self-convolutive recurrence}

\begin{abstract}

We consider a self-convolutive recurrence whose solution is the sequence of coefficients in the asymptotic expansion of the logarithmic derivative of the confluent hypergeometic function $U(a,b,z)$. By application of the Hilbert transform we convert this expression into an explicit, non-recursive  solution in which the $n$th  coefficient is expressed as the $(n-1)$th moment of a measure, and also as the trace of the $(n-1)$th iterate of a linear operator. Applications of these sequences, and hence of the explicit solution provided, are found in quantum field theory as the number of Feynman diagrams of a certain type and order, in Brownian motion theory, and in combinatorics.
\newline
\emph{Aequationes Mathematicae}, 80, 291--318 (2010).

\end{abstract}

\maketitle

\section{Introduction}

We study the sequence $\sequence(\alpha_1,\alpha_2,\alpha_3)=(u_n)_{n=1}^\infty$ defined by the self-convolutive recurrence
\begin{equation}
u_n = (\alpha_1 n + \alpha_2) u_{n-1} + \alpha_3 \sum_{j=1}^{n-1} u_j u_{n-j}, \qquad u_1=1,
\label{eq:u_n}
\end{equation}
and show how to derive a closed-form solution as a Mellin transform, i.e.\ one of the form
\begin{equation}
u_n = \int_0^\infty x^{n-1} \mu(x) \,dx,
\label{eq:Mellin}
\end{equation}
at least, in a reasonably broad special case (in general, there are extra terms, but we give an explicit expression for them).
In other words,  $u_n$ is the $(n-1)$th moment of a measure.
Clearly a representation such as (\ref{eq:Mellin}) is useful if one wishes to compute $u_n$ for large $n$---as it is nonrecursive---or study the asymptotics via Laplace's method \cite{Bender78}. Incidentally the representation additionally defines an analytic continuation of $u_n$ to $\Real n > 1$. 

Our study of such sequences arose from a particular example, namely the coefficients $(u_n)$ of the asymptotic expansion of the logarithmic derivative of the Airy function:
\[
\frac{\Ai'(z)}{\Ai(z)} \sim 2\sqrt{z} \sum_{n=0}^\infty \frac{(-)^nu_n}{(4z)^{3n/2}}
\]
with $u_0=-\half$. As $\Ai'/\Ai$ obeys a Riccati equation, the $(u_n)$ form a recurrence of the type given in (\ref{eq:u_n}):
\[
u_n = (6n-8)u_{n-1} + \sum_{j=1}^{n-1} u_j u_{n-j}.
\]
It was established in \cite{Kearney07}, thereby settling a long-standing question, that
\[
u_n = \frac{1}{\pi^2} \int_0^\infty \frac{x^{n-4/3}}{\Ai(\frac{1}{4}x^{2/3})^2 + \Bi(\frac{1}{4}x^{2/3})^2 } \,dx, \quad n\ge 1,
\]
which is of the form (\ref{eq:Mellin}) with $\mu(x)=\pi^{-2}x^{-1/3}\big/\big[\Ai(\frac{1}{4}x^{2/3})^2 + \Bi(\frac{1}{4}x^{2/3})^2\big]$.
This sequence and a related one described later are encountered in a number of applications including the moments of the area function of a Brownian excursion \cite{Kearney07,Kearney09}, quantum field theory \cite{Cvitanovic78} (quite a fertile ground for these: more of this shortly) and combinatorics \cite{Janson03}.
For reference this sequence is \#A062980 in the Online Encyclopedia of Integer Sequences (OEIS) \cite{Sloane73}, and begins
\[
S(6,-8,1) = (1,5,60,1105,27120,828250,\ldots).
\]

An obvious question to ask was whether this result can be generalised, and the answer is yes. A trivial example is $\sequence(1,\alpha,0)$ for which the $n$th term is $\Gamma(n+\alpha+1)/\Gamma(\alpha+2)$ and in Mellin form is 
\[
\sequence(1,\alpha,0)_n = \int_0^\infty x^{n-1} \frac{x^{\alpha+1}e^{-x} }{\Gamma(2+\alpha)} \,dx.
\]
However, it has a more interesting `relative' (we shall say later in what sense)
\[
\sequence(1,-2,1)_n= \int_0^\infty \frac{x^{n-2} e^x}{\Ei(x)^2 + \pi^2 } \,dx,
\]
where $\Ei$ denotes the exponential integral function \cite{Abramowitz64}. It turns out that this sequence has been studied extensively in combinatorics and algebra \cite{Flajolet09}, because it is the number of \emph{connected}, or \emph{indecomposable}, permutations of $[1\ldots n]$, i.e.\ those not fixing $[1\ldots j]$ for $0<j<n$. For reference this sequence is \#A003319 in the OEIS and begins
\[
S(1,-2,1) = (1,1,3,13,71,461,3447,\ldots).
\]
The integral representation is, we think, a new result. 

There are more known cases in this category, in fact, and here is another result from combinatorics:
\[
\sequence(2,-3,1)_n= \sqrt{\frac{2}{\pi}} \int_0^\infty \frac{ x^{n-3/2}  e^{x/2}}{G_1(\half x)^2 + \pi } \,dx,
\]
where
\[
G_1(x) \equiv   \sum_{r=0}^\infty \frac{x^{r+\half}}{(r+\half)r!} .
\]
This sequence (\#A000698 in the OEIS) begins
\[
S(2,-3,1) = (1, 2, 10, 74, 706, 8162, 110410, \ldots)
\]
and is the number of nonisomorphic connected Feynman diagrams of order $2(n+1)$ arising in a simplified model of quantum electrodynamics (QED) \cite{Cvitanovic78}. On multiplication term-by-term by $2n-1$ the following sequence (\#A005416) is generated,
\[
(1, 6, 50, 518, 6354, 89782, \ldots)
\]
which is the number of `vertex graphs' of order $2n$ arising in the QED perturbation series for the electron  magnetic moment \cite{Cvitanovic77}.
The above expression as the moment of a distribution has been stated by R.~Groux \cite{Sloane73} though the result and derivation do not seem to have been formally published.
This is one of the very few cases that have been solved.

A sequence that bears intriguing similarity to that one is
\[
\sequence(2,-2,1)_n=  \int_0^\infty \frac{ 2x^{n-1}  e^{x/2}}{G_2(\half x)^2 + \pi^2 G_3(\half x)^2} \,dx,
\]
where
\[
G_2(x) \equiv  2 - G_3(x)\ln x +   \sum_{r=0}^\infty  \frac{C_r d_r x^{r+1}}{r!4^r}
 , \quad
G_3(x) \equiv   \sum_{r=0}^\infty  \frac{C_r x^{r+1}}{r!4^r}  \,;
\]
here $C_r=\frac{1}{r+1} {2r \choose r}$  is the $r$th Catalan number and 
\[
d_r = \textstyle  -\gamma+2\ln2+ \frac{1}{r+1} + 2\sum_{j=1}^r \left(\frac{1}{j}-\frac{1}{2j-1}\right)  
\]
where $\gamma=0.5772\ldots$ is Euler's constant.
This sequence (\#A005412 in the OEIS) begins
\[
S(2,-2,1) = (1, 3, 18, 153, 1638, 20898, \ldots)
\]
and is the number of Feynman diagrams with proper self-energies arising in QED \cite{Cvitanovic78}; the integral representation is apparently new. 

Another example which generates a neater result is this, invoking the modified Bessel functions. The conciseness of its form is entirely due to the familiarity of the functions involved, though, and in fact it is almost identical in construction to the previous one:
\begin{eqnarray*}
S(2,-4,1)_n &=&  \int_0^\infty  \frac{2x^{n-2} \, dx}{K_0(\frac{1}{4} x)^2 + \pi^2 I_0(\frac{1}{4} x)^2} \\
S(2,0,-1)_n &=&  \int_0^\infty  \frac{2x^{n-2} \, dx}{K_1(\frac{1}{4} x)^2 + \pi^2 I_1(\frac{1}{4} x)^2} .
\end{eqnarray*}
The first few terms are
\begin{eqnarray*}
S(2,-4,1) &=& (1,1,4,25,208,2146,26368,\ldots) \\
S(2,0,-1) &=& (1,3,12,63,432,3798,41472,\ldots).
\end{eqnarray*}
The first of these occurs in QED as the number of Feynman diagrams with exact propagators \cite{Cvitanovic78} and is \#A005411 in the OEIS, though the expression as an integral is, again, apparently new.
We are not sure if the second expression, which is not in the OEIS at the time of writing, has an application in the same field.


This paper is organised as follows: the next section gives a complete exposition of the methods, and the one after that is devoted to particular examples. Notation is standard throughout; $\N$ denotes the set of natural numbers $1,2,\ldots$; $-\N$ denotes their negatives; $\N_0$,$-\N_0$ denote the same sets including zero. We use $\digamma(x)$ to denote the digamma function, $\Gamma'(x)/\Gamma(x)$; note that for $m\in\N_0$, we have $\digamma(1+m) = -\gamma + \sum_{r=1}^{m} \frac{1}{r}$ and $\digamma(\half+m)=-\gamma-2\ln2+\sum_{r=1}^m \frac{2}{2r-1}$, where $\gamma=0.5772\ldots$ is Euler's constant and the sums are disregarded for $m=0$.

\section{Methods}


\subsection{Connection with Riccati equation, and solution as integral}

The use of the Riccati equation is fundamental, and we argue as follows. The sequence is generated from the asymptotic expansion of a function $\myfunc$, i.e.\
\begin{equation}
\myfunc(z) \sim \sum_{n=0}^\infty \frac{(-)^n u_n}{z^n}, \qquad |\arg z| < \pi,
\label{eq:asy}
\end{equation}
where we have added in the 0th term defined in the obvious way as $u_0=\frac{1}{\alpha_1+\alpha_2}$ (if $\alpha_1+\alpha_2=0$, we consider the limit $\alpha_1+\alpha_2\to0$). Note however that when subsequently we give formulas for $u_n$, they will not necessarily be valid for $n=0$. The function $\myfunc$ is assumed to be single-valued in the complex plane cut along the negative real axis, and in the situation at hand this condition will be satisfied in all cases bar one (the `algebraic case'), which we deal with separately.

The recursion in the $(u_n)$ corresponds to the following Riccati equation $\mathcal{R}(\beta_1,\beta_2,\beta_3,\beta_4)$ defined below (\ref{eq:ricc}), in the sense that if we find a function satisfying (\ref{eq:ricc}) that has an asymptotic
expansion (\ref{eq:asy}), then the coefficients will satisfy (\ref{eq:u_n}) and thus be
the sought sequence:
\begin{equation}
\myfunc'(z) = \beta_1 \myfunc(z)^2 + (\beta_2+\beta_3/z) \myfunc(z) + \beta_4
\label{eq:ricc}
\end{equation}
with coefficients related by
\begin{equation}
\beta_1=-\frac{\alpha_3}{\alpha_1}, \quad \beta_2=\frac{\alpha_1+\alpha_2+2\alpha_3}{\alpha_1(\alpha_1+\alpha_2)}, \quad \beta_3=1+\frac{\alpha_2}{\alpha_1}, \quad \beta_4 = -\frac{\alpha_1+\alpha_2+\alpha_3}{\alpha_1(\alpha_1+\alpha_2)^2}
\end{equation}
(again, the case $\alpha_1+\alpha_2=0$ can be dealt with by a limiting argument).
Suppose we have found a solution to this equation with the appropriate regularity properties (in particular, $\lim_{x\to\infty}\myfunc(x)= \frac{1}{\alpha_1+\alpha_2}$).
Then we only have to extract the $(u_n)$ by the Fourier integral:
\[
(-)^n u_n = \frac{1}{2\pi\I} \int_\mathcal{C} \myfunc(z) z^n \, \frac{dz}{z}
\]
where the contour $\mathcal{C}$ runs anticlockwise in a large circle from $-\infty-\I\varepsilon$ to $-\infty+\I\varepsilon$ (but it is not closed).

It is probably fair to say that this last step has been something of a stumbling-block to explorers, and is the reason that the results in this paper have lain undiscovered for so long. On the face of it, the right-hand side of (\ref{eq:asy}) is not a meaningful function, so apparently it cannot be valid to multiply by $z^{n-1}$ and then perform the integral in such a cavalier fashion. However, the asymptotic series may be written as a finite sum of $N$ terms plus a remainder, $R_N(z)$ say, and crucially this remainder is bounded by  $c_N|z|^{-N-1}$ in $\{|z|>\rho,\;|\arg z|<\pi\}$ for some constant $c_N$ independent of $z$. The integral can then be performed around a contour of finite radius $>\rho$ and the contour expanded to infinite radius, causing the remainder term to vanish. Finally $N$ is arbitrary and therefore can be made as big as we please. Incidentally we consider that much of the difficulty is caused by the `$\sim$' notation, which makes asymptotic analysis look like a branch of witchcraft not possessing the rigour of Taylor series. In fact, there is a symmetry between the two notions: whereas an asymptotic series has a remainder term that is bounded for a given $N$ for all $z$ inside some annulus, the Taylor series has a remainder term that is bounded for a given $z$ for all $n$ exceeding some $N$. For our purposes it is the first one that is useful.
It is worth emphasizing that we are assuming the asymptotic expansion (\ref{eq:asy}) to hold uniformly in the whole cut plane. This is not a problem in the case considered here, but for extensions or generalisations of the methods shown the region of validity needs to be borne in mind.

 This contour integral can now be evaluated as the sum of three parts, as follows. We collapse the contour on to the negative real axis: the resulting integral is the first part. Secondly, in wrapping the contour round the branch cut $\R^-$ we may generate a contribution from the origin if the integrand is not integrable there. Finally, along the way we may have picked up contributions from other singularities of $\myfunc$ anywhere in the complex plane except the negative real axis and the origin.

It is convenient to classify cases by reference to their singularities.

\begin{defn} \label{def:sr}
The sequence (or, equivalently, the associated Riccati equation) is said to be \emph{simply represented} if it obeys both the following conditions: 
\begin{itemize}
\item[(i)] $\myfunc$ has no singularities other than simple poles in the cut plane $\C \setminus (-\infty,0]$; 
\item[(ii)] $\lim_{z\to0}z\myfunc(z)$ exists and is finite.
\end{itemize}
\end{defn}

In the context of this paper there are very few sequences that are not simply represented, and we can and do deal with them individually. We subdivide the simply represented sequences as follows:

\begin{defn}
A simply represented  sequence (or Riccati equation)  is said to be \emph{regular} if it obeys both the following stronger conditions: 
\begin{itemize}
\item[(i*)] $\myfunc$ has no singularities in the cut plane $\C \setminus (-\infty,0]$;
\item[(ii*)] $\lim_{z\to0} z\myfunc(z)=0$.
\end{itemize}
It is said to be \emph{quasiregular} if it obeys (ii) but not (ii*). Otherwise 
it is said to be \emph{irregular}.
\end{defn}

Incidentally all the examples mentioned in the Introduction are regular. In the regular case we can collapse the contour around the branch cut and write $z=-x$ to give
\begin{equation}
u_n = \frac{1}{2\pi\I} \int_{\to 0}^\infty \diff{\myfunc}(x) x^{n-1} \, dx 
\label{eq:regsoln}
\end{equation}
where $\diff{\myfunc}(x)$ denotes the jump in $\myfunc(z)$ across the branch cut at $z=-x$, i.e.\  $\lim_{\varepsilon\to0+} [\myfunc(-x+\I\varepsilon) - \myfunc(-x-\I\varepsilon)]$.
So for a regular sequence, we have now obtained a solution in the form (\ref{eq:Mellin}) with
\begin{equation}
\mu(x) = \frac{1}{2\pi\I} \diff{\myfunc}(x). 
\end{equation}
Note that $\myfunc$ can be recovered from $\diff{\myfunc}$ immediately through the Hilbert transform:
\[
\myfunc(z) = -\frac{1}{2\pi\I} \int_0^\infty \frac{\diff{\myfunc}(x)}{x+z} \,dx, \qquad \Real z>0,
\]
with the definition elsewhere being provided by analytic continuation, paying attention to the need to cut the plane along the negative real axis.

In the quasiregular case we have an extra term:
\[
u_n = \frac{1}{2\pi\I} \int_{\to 0}^\infty \diff{\myfunc}(x) x^{n-1} \, dx 
+ \left\{ \begin{array}{rl} -r^\circ, & n=1 \\ 0, & n>1 \end{array}\right\}
\]
where $r^\circ$ is the residue at the origin.
If desired this can be merged into $\mu(x)$ as a delta-function, $-r^\circ\cdot\delta(x)$.

If there are simple poles in the cut plane, located at $z=\zeta_j$ and of residue $s_j$, further terms are added:
\begin{equation}
u_n = \frac{1}{2\pi\I} \int_{\to 0}^\infty \diff{\myfunc}(x) x^{n-1} \, dx  -r^\circ\cdot0^{n-1}
 - \sum_j s_j (-\zeta_j)^{n-1}.
\label{eq:gensoln}
\end{equation}
Writing the contribution from the origin as $-r^\circ\cdot 0^{n-1}$, interpreting $0^0=1$,  makes it look like the contributions arising from other poles. 


Incidentally, simple poles \emph{on} the negative real axis make their presence felt when the integral is collapsed on to the branch cut. 
A direct calculation of the expression $\lim_{\varepsilon\to0+} [\myfunc(z+c+\I\varepsilon) -\myfunc(z+c-\I\varepsilon)]$ shows it to be equal to $-2\pi\I\delta(x-c)$, with $x=-z\in\R^+$, but a quicker route is to use the Hilbert transform:
\[
\int_0^\infty \frac{\delta(x-c)}{z+x} \,dx = \frac{1}{z+c}
\]
showing that the right interpretation of a pole in $\R^-$ in the generating function $\myfunc$ is a delta-function in the measure $\mu$, of strength equal to minus the residue of the pole.

If we were to relax the requirement (ii) in Definition \ref{def:sr}, to allow a singularity worse than a simple pole, but still of power order, at the origin, then we would have to deal with a contribution from there of the form
\[
\frac{1}{2\pi\I}\int_{\varepsilon e^{-\pi\I}}^{\varepsilon e^{\pi\I}} \myfunc(z) z^n \, \frac{dz}{z}
\]
i.e.\ the integral around a small loop. By stipulating that $z^m\myfunc(z)$ be bounded at the origin for sufficiently large $m$, we ensure that a singularity of this type can affect only finitely many of the $(u_n)$. We then say that the sequence is \emph{almost simply represented}.

We have shown in general how to solve for $u_n$, but now need to give an explicit expression for $\myfunc$, so we attend to that next.

\subsection{Solving the Riccati equation}

First, we exclude the case $\alpha_1=0$, as then the Riccati equation is purely algebraic.
 This is better treated on its own, which we do in \S\ref{sec:Alg}.

The next special case that requires attention is $\alpha_3=0$. In that case the Riccati equation is a linear differential equation and $\myfunc$ is essentially an exponential integral function \cite[\S5]{Abramowitz64}. We deal with this as the first of the special cases in the next section. It is well known that this function is regular in $\C\setminus(-\infty,0]$ and has a singularity of power order at worst at the origin.

Another observation is the `scaling law'
\[
S(t\alpha_1,t\alpha_2,t\alpha_3)_n = t^{n-1} S(\alpha_1,\alpha_2,\alpha_3)_n.
\]
By solving for $(-)^{n-1}u_n$ instead of $u_n$, we can assume that $\alpha_1>0$.

Having dealt with these cases we are now free to assume $\beta_1\ne0$, and so we can write $\myfunc(z)=-w'(z)/\beta_1 w(z)$ in the usual way to obtain the linear differential equation
\begin{equation}
w'' - (\beta_2+\beta_3/z) w' + \beta_1 \beta_4 w = 0.
\end{equation}
The solution of this is essentially a confluent hypergeometric function (CHGF; \cite[\S13]{Abramowitz64}), and the appropriate solution is
\[
w(z) \propto e^{-(ak/b)z}U(a,b,kz)
\]
where
\begin{equation}
k 
 = \frac{1}{\alpha_1}, \quad a = -\beta_1 = \frac{\alpha_3}{\alpha_1},   \quad b=-\beta_3=-1-\frac{\alpha_2}{\alpha_1};
\end{equation}
notice that by assumption $k,a,b$ are finite, $k>0$ and $a\ne 0$. The effect of `scaling' (q.v.) on these parameters is $k\mapsto t^{-1}k$, with $a,b$ fixed.
Note that there is an apparent difficulty when $b=0$, as then $e^{-(ak/b)z}$ is then undefined. However, the effect of the exponential is simply to add a constant on to $w'/w$, and in the asymptotic expansion this only affects $u_0$ which is formally infinite, but not of interest. Unlike some of the other special cases (notably $\alpha_1=0$, which requires special treatment), this can be dealt with by simply letting $b\to0$ in the final expression.
The function $U$ is given (for $b\notin\Z$) by
\begin{equation}
U(a,b,z) = \frac{\Gamma(1-b)}{\Gamma(a-b+1)} M(a,b,z) + \frac{\Gamma(b-1)}{\Gamma(a)} z^{1-b}M(a-b+1,2-b,z)
\label{eq:UfromM}
\end{equation}
in which $M(a,b,z)$, also often denoted $_1F_1(a;b;z)$, is the Kummer function, given by the Taylor series $\sum_{r=0}^\infty \frac{(a)_r}{(b)_r}\frac{z^r}{r!}$, and is the regular solution of the ODE
\[
zM''+(b-z)M'-aM=0.
\]
Note that
\[
U(a,b,z) = z^{1-b} U(a-b+1,2-b,z),
\]
 which is one of `Kummer's transformations'.
A contour integral representation of $U(a,b,z)$ for $\Real z>0$ is
\[
U(a,b,z)=\frac{1}{\Gamma(a)}\int_0^\infty e^{-zt} t^{a-1}(1+t)^{b-a-1}\,dt\qquad (\Real a > 0)
\]
or
\[
U(a,b,z)=\frac{-\Gamma(1-a)e^{-\pi\I a}}{2\pi\I} \int_{\infty}^{(0-)} e^{-zt} t^{a-1}(1+t)^{b-a-1}\,dt\qquad (a \notin \N)
\]
where in the second expression the plane is cut from 0 to $+\infty$ and $-1$ to $-\infty$ and the contour runs from $+\infty$ just below the cut, leftwards, round the origin and then back to $+\infty$ on the upper side of the cut.
The asymptotic expansion of $U$ is
\begin{equation}
U(a,b,z) \sim z^{-a} \sum_{n=0}^\infty \frac{(-)^n (a)_n(a-b+1)_n}{n!\,z^n}, \qquad |\arg z| < \textstyle \frac{3}{2} \pi,
\label{eq:asexpU}
\end{equation}
so the behaviour of $U$ is known around the contour $\mathcal{C}$ mentioned earlier.

We now need to evaluate $\diff{\myfunc}(x)$. From $\myfunc(z)=-w'(z)/\beta_1 w(z)$ we deduce
\[
\myfunc(z) = \frac{k}{a}  \frac{\frac{d}{d(kz)}U(a,b,kz)}{U(a,b,kz)} - \frac{k}b.
\]
Let $\Uplus(a,b,z)$ and $\Uminus(a,b,z)$ denote the values of the $U(a,b,z)$ for $\arg(z)=+\pi$ and $-\pi$ respectively. Then for $z\in\R^-$,
\[
\myfunc^\oplus(z)-\myfunc^\ominus(z) = -\frac{k}{a} \left. \frac{\Wr{\Uplus}{\Uminus}}{\Uplus\Uminus} \right|_{kz}
\]
with $\Wr{}{}$ denoting the Wronskian. Now
\[
\displaystyle U^\oplus (a,b,z) = 
 \frac{\Gamma(1-b)}{\Gamma(a-b+1)} M(a,b,z) + \frac{\Gamma(b-1)}{\Gamma(a)} (z^{1-b})^\oplus M(a-b+1,2-b,z) 
\]
and similarly for the lower branch $U^\ominus$. So the Wronskian of $U^\oplus$ and $U^\ominus$ can be dealt with easily enough using the identity
\[
\Wr{M(a,b,z)}{z^{1-b}M(a-b+1,2-b,z)} = (1-b)z^{-b}e^z,
\]
and the result emerges as, on writing $x$ for $-z$,
\[
\Wr{\Uplus(a,b,z)}{\Uminus(a,b,z)} = \frac{-2\pi\I x^{-b}e^{-x}}{\Gamma(a)\Gamma(a-b+1)}.
\]
Hence
\begin{equation}
\mu(x)=\frac{1}{2\pi\I} \diff{\myfunc}(x) = \frac{ k (kx)^{-b}e^{-kx}}{\Gamma(a+1)\Gamma(a-b+1)} \cdot \frac{1}{\Uplus(a,b,-kx)\Uminus(a,b,-kx)}.
\label{eq:result1}
\end{equation}
A little more work is needed on the last part. This involves the analytic continuation of $U$ to the left half-plane, which is:
\begin{eqnarray}
U(a,b,xe^{\pm\pi\I}) &=& e^{-x} \left[\frac{\Gamma(1-b)}{\Gamma(a-b+1)} M(b-a,b,x) \right.   \nonumber
\\ && \;\;\;\;\;\;\; \left. - \, e^{\mp\pi\I b}x^{1-b} \frac{\Gamma(b-1)}{\Gamma(a)} M(1-a,2-b,x)\right]
\end{eqnarray}
We can therefore write
\[
U(a,b,xe^{\pm \pi\I}) = U_R(x) \pm \I U_I(x)
\]
with
\begin{eqnarray}
U_R(x) &=& e^{-x} \left[\frac{\Gamma(1-b)}{\Gamma(a-b+1)} M(b-a,b,x) \right. \nonumber \\
&& \;\;\;\;\;\;\; \left. -\, (\cos \pi b )\, \frac{\Gamma(b-1)}{\Gamma(a)} x^{1-b} M(1-a,2-b,x)\right] \nonumber \\
U_I(x) &=& (\sin \pi b)\,  \frac{\Gamma(b-1)}{\Gamma(a)} e^{-x}   x^{1-b} M(1-a,2-b,x)
\label{eq:URI}
\end{eqnarray}
and the $\Uplus\Uminus$ term from above is expressed as
\begin{equation}
\Uplus(a,b,-kx)\Uminus(a,b,-kx) = U_R(kx)^2+U_I(kx)^2.
\label{eq:URI2}
\end{equation}

Now that we have solved for $\myfunc$ we can be rigorous about its singularities.

\begin{lem} \label{lem:cp}
Suppose $\alpha_1\ne0$. The only singularities of $\myfunc$ in the cut plane, if any, are simple poles of residue $\alpha_1/\alpha_3$; unless $\alpha_3=0$, in which case there are none.
\end{lem}
\noindent Proof.
Other than at the origin, the CHGF can have   no singularities and nor can it have zeros of order $>1$, as either would  violate the differential equation that it satisfies. Hence the only singularities that $w'/w$ can have, except at the origin, are simple poles of residue 1. So the residues of $\myfunc$ are always $-1/\beta_1$, unless $\beta_1=0$ in which case $\myfunc$ is an exponential integral and is free from singularities except possibly at the origin. $\Box$

Now let us be more precise about the singularity of $\myfunc$ at the origin. By examining $U(a,b,z)$ there we can deduce:
\begin{lem} \label{lem:origin}
Notation as above. Writing $r^\circ=\lim_{z\to0}z\myfunc(z)$, we have:
\begin{itemize}
\item Case $a\notin-\N_0$, $a-b+1\notin-\N_0$. If $b>1$ then $r^\circ=(1-b)/a$; otherwise, $r^\circ=0$. 
\item Case $a\in-\N_0$, $a-b+1\notin-\N_0$. Then $r^\circ=0$.
\item Case $a\notin-\N_0$, $a-b+1\in-\N_0$. Then $r^\circ=(1-b)/a$.
\item Case $a\in-\N_0$, $a-b+1\in-\N_0$ (so $b\in\Z$). If $b\le 1$ then $r^\circ=(1-b)/a$; otherwise $r^\circ=0$. 
\end{itemize}
\end{lem}
\noindent Proof. 
The first case can be obtained from (\ref{eq:UfromM}), examining which of the two halves predominates when $z\to 0$; see also \cite[\S13]{Abramowitz64}. For $a\in-\N_0$,  $U(a,b,z)$ is a polynomial.
The second and fourth cases follow from direct examination of the coefficients. The third follows from Kummer's transformation to exchange $a$ for $a-b+1$, picking up a factor of $z^{1-b}$ in the process. $\Box$

Incidentally if we assume that $\myfunc(z)$ looks like $c/z$ at the origin then we can deduce immediately that $c=0$ or $c=(1-b)/a$ directly from the Riccati equation; but choosing the right one is not quite so obvious.

We have arrived at:
\begin{thm}
If $\alpha_1>0$ and $\alpha_3\ne0$ then  $S(\alpha_1,\alpha_2,\alpha_3)$ is simply represented, and hence solved by $(\ref{eq:gensoln})$. $\Box$
\end{thm}

Somewhat paradoxically, the ostensibly trivial case $\alpha_1=0$ has pecularities of its own that do not always permit a simple representation. We leave this case until the end of the paper because it bears resemblances to the general case that will only become apparent when a good selection of special cases have been exposed.

Now we have to compute all the various bits.
The following result summarises what we have found about $\mu(x)$ and shows it to be a measure, up to a prefactor (remember that $k>0$ by assumption):

\begin{thm}
The function $\mu$ in $(\ref{eq:Mellin})$ is given by $(\ref{eq:result1},\ref{eq:URI},\ref{eq:URI2})$. After extraction of a constant factor
it is positive for all $x>0$.
$\Box$
\end{thm}

Now we turn to the question of zeros of $U$ in the cut plane, which give rise to simple poles in $\myfunc$.
Denote by $\Delta$ the change in argument, divided by $2\pi$, of $U(a,b,z)$ as $z$ passes along each of the following paths:  $\Delta_\mathcal{C}$, around the loop $\mathcal C$ anticlockwise from $-\infty$ to $-\infty$ as previously defined; $\Delta_+$, along the path from  $-\infty$ to $0$ above the branch cut; $\Delta_0$, for the clockwise loop around the origin; $\Delta_-$, along the path from $0$ to $-\infty$ below the cut. 
Then by the Argument Principle, the number of zeros of $U$ in the cut plane, counted according to multiplicity, is
\[
\# U = \Delta_\mathcal{C} + \Delta_{+} + \Delta_0 + \Delta_{-}
\]
(remember that $U$ has no poles in there). As the properties of $U(a,b,z)$ on $\mathcal C$ and near the origin are well documented \cite[\S 13]{Abramowitz64}, we can track $\arg U$ very easily provided we can find $\Delta_+$ and $\Delta_-$, which we expect to be able to do from examination of $U_R(x)$ and $U_I(x)$, in certain cases at least. Note that $\Delta_+=\Delta_-$ by symmetry.

\begin{lem} \label{lem:zc1}
The number of zeros of $U(a,b,z)$ in the cut plane is as follows.
\begin{itemize}
\item[(i)] If $-2m-1 \le a<-2m+1$, $m\in\N_0$, and $b<a+1$, then $\# U=2m$.
\item[(ii)] If $a=b\ge 1$ then $\# U=0$.
\end{itemize}
\end{lem}
\noindent Proof. In all cases one has $\Delta_\mathcal{C}=-a$, by (\ref{eq:asexpU}).

\notthis {
\noindent (i) The assumptions imply positivity of $M(1-a,2-b,x)$ for $x>0$, so
\[
U_I(x) = \frac{-\pi}{\Gamma(2-b)\Gamma(a)} x^{1-b} e^{-x} M(1-a,2-b,x) <0.
\]
As $z$ passes along $\mathcal C$, $\arg U$ changes from $+\pi a$ to $-\pi a$ (and $U$ moves clockwise). Above the branch cut, $\Imag U<0$. At the origin, $U$ tends to the \emph{positive} real value $\Gamma(1-b)/\Gamma(a-b+1)$. Below the branch cut, $\Imag U>0$. So in total, $\# U = -a + \half a + 0 + \half a =0$. 

\noindent(ii)
Again the assumptions imply positivity of $M(1-a,2-b,x)$ for $x>0$, so
\[
U_I(x) = \frac{-\pi}{\Gamma(2-b)\Gamma(a)} x^{1-b} e^{-x} M(1-a,2-b,x) >0.
\]
As $z$ passes along $\mathcal C$, $\arg U$ changes from $+\pi a$ to $-\pi a$ (so $U$ moves \emph{anticlockwise}, as $a<0$). Above the branch cut, $\Imag U>0$. At the origin, $U$ tends to the \emph{positive} real value $\Gamma(1-b)/\Gamma(a-b+1)$. Below the branch cut, $\Imag U<0$. So again $\# U = -a + \half a + 0 + \half a =0$.

} 

\noindent(i) We deal with the case $b<1$ first.
Write $a=-2m-1+\delta$, with $0<\delta<2$, thereby temporarily excluding the case when $a$ is a negative odd integer.
 The assumptions imply positivity of $M(1-a,2-b,x)$ for $x>0$, so
\[
U_I(x) = \frac{-\pi}{\Gamma(2-b)\Gamma(a)} x^{1-b} e^{-x} M(1-a,2-b,x) \left\{ \begin{array}{rr} <0, & 0<\delta< 1\\ >0, & 1<\delta<2 \end{array} \right.
\]
Note that $U(0)=\Gamma(1-b)/\Gamma(a-b+1)>0$.
From (\ref{eq:asexpU}), and having just deduced that $\Imag U(a,b,z)$ does not change sign as $z$ passes from $\infty e^{+\pi\I}$ to 0 above the cut, we must have $\Delta_+=\half(\delta-1)$. So $\# U = -a + \half(\delta-1) + 0 + \half (\delta-1) =2m$.
 The case when $a$ is a negative odd integer is more subtle because the number of zeros changes by two on account of a pair of zeros approaching the negative real axis. As the complement of the cut plane (in the Riemann surface) is a closed set, it contains its limit points, and so the extra pair of zeros remains in that set as $a$ approaches and touches a negative odd integer from below. Hence, for example, with $a=-3$ we have $\# U =2$ and not 4.
(When $a$ is a negative even integer, i.e.\ $\delta=1$, $U(z)$ simply travels along the positive real axis without going through the origin, so there is no difficulty.)

The remaining case to deal with is $0<a<1$,  $1<b<a+1$. Then $U_I$ is negative above the branch-cut, by the same reasoning as above. The behaviour at the origin is different, as now $U(a,b,z)\sim \frac{\Gamma(b-1)}{\Gamma(a)} z^{1-b}$; this implies $\Delta_0=b-1$, and it is easily seen that $\Delta_+=\half(a-b+1)$. Thus $\# U=0$ again, just as we found above for $0<a<1$ and $b<1$.

\noindent(ii)
From the behaviour of $U$ at $0$, we have $\Delta_0=a-1$. To attend to $\Delta_+$ and $\Delta_-$ we consider the behaviour of $U(a,a,z)e^{-\pi \I(1-a)}$, first on the upper side of the cut. This goes from $-0$ to $+\infty$ as $z$ runs from $\infty e^{+\pi\I}$ to $0 e^{+\pi\I}$, and we need to ascertain which side of the origin it passes on the way. By the analytic continuation formula for $U$ we have, again writing $z=-x$, $x\in\R^+$,
\[
U(a,a,x e^{+\pi\I})e^{-\pi \I(1-a)} = -e^{-x}\Gamma(1-a)e^{\pi\I a}M(0,1,x) - \mbox{(real)},
\]
so its imaginary part is $-\pi e^{-x}/\Gamma(a)<0$. Accordingly, $\Delta_+=+\half$. So $\# U= -a + \half + (a-1) + \half =0$.
$\Box$

\subsection{Symmetries}

Two symmetries are worth mentioning. First, by Kummer's transformation we have
\[
\frac{1}{a} \frac{U'(a,b,z)}{U(a,b,z)} = \frac{a-b+1}{a} \left( \frac{1}{a-b+1} \frac{U'(a-b+1,2-b,z)}{U(a-b+1,2-b,z)}\right) + \frac{1-b}{az}.
\]
Accordingly
\begin{eqnarray}
S(\alpha_1,\alpha_2,\alpha_3)_n &=& \textstyle \frac{2\alpha_1+\alpha_2+\alpha_3}{\alpha_3}  \, S(\alpha_1,-4\alpha_1-\alpha_2,2\alpha_1+\alpha_2+\alpha_3)_n  \nonumber \\
&& - 
 \left\{ \begin{array}{rl} \frac{2\alpha_1+\alpha_2}{\alpha_3}, & n=1 \\ 0, & n>1 \end{array}\right\}
\end{eqnarray}
which we call the Kummer symmetry.
Two sequences can thereby be solved `for the price of one'; as an example,
\[
S(6,-8,1)_n = 5\, S(6,-16,5)_n - 4 \cdot 0^{n-1}.
\]

Secondly, the set of Riccati equations $\mathcal{R}(\ldots)$ as previously defined is closed under reciprocation, i.e.\ the transformation $\myfunc(z)\mapsto \theta/\myfunc(z)$, which effects the following transformation,
\[
\mathcal{R}(\beta_1,\beta_2,\beta_3,\beta_4) \mapsto 
\mathcal{R}(-\beta_4/\theta,-\beta_2,-\beta_3,-\beta_1\theta);
\]
here $\theta = -u_0^2=-(\alpha_1+\alpha_2)^{-2}$ to make sure that $u_1=1$ in the transformed sequence.
This induces the following transformation on the space of sequences:
\[
\sequence(\alpha_1,\alpha_2,\alpha_3) \mapsto
\sequence(\alpha_1,-2\alpha_1-\alpha_2,\alpha_1+\alpha_2+\alpha_3).
\]
Writing $^\dagger$ for the reciprocal, we have for example
\begin{eqnarray*}
S(1,\alpha-1,0)^\dagger &=& S(1,-1-\alpha,\alpha)\\
S(6,-8,1)^\dagger &=& S(6,-4,-1).
\end{eqnarray*}
The effect on the CHGF parameters is
\[
k^\dagger = k ; \quad a^\dagger = a-b ; \quad b^\dagger = -b,
\]
and so $2a^\dagger-b^\dagger=2a-b$; also sequence is self-reciprocal iff $b=0$. Clearly $(S^\dagger)^\dagger \equiv S$. Note carefully that the reciprocal of a regular sequence need not be regular.

\subsection{Asymptotics}

The expression as the $n$th moment of a measure makes it particularly easy to examine asymptotics.
For example, if 
\[
\mu(x)\sim Ax^\nu e^{-\lambda x} \quad (x\to\infty)
\]
then
\[
\int_0^\infty x^{n-1} \mu(x)\,dx \sim A\Gamma(n+\nu) \lambda^{-n-\nu}.
\]
From (\ref{eq:asexpU},\ref{eq:result1}) we have
\[
\mu(x) \sim \frac{k(kx)^{2a-b}e^{-kx}}{\Gamma(a+1)\Gamma(a-b+1)}, \qquad x\to\infty,
\]
provided the Gamma terms do not vanish, allowing us to write down:
\begin{thm}\label{thm:asu}
If $S$ is simply represented, $a\notin-\N$, $a-b\notin-\N$, $k\ne 0$, then
\[
u_n\sim  \frac{\Gamma(n+2a-b)}{\Gamma(a+1)\Gamma(a-b+1)k^{n-1}}  , \qquad n\to\infty.
\]
Further, the reciprocal sequence to $S$ has the same behaviour at leading order.
\end{thm}
\noindent Proof. Clear from the above. Any poles in $\myfunc$ generate power order behaviour, which is asymptotically dominated by the Gamma function. $\Box$ 

When $a\in-\N$ the function $U(a,b,z)$ is a polynomial of degree $-a$ and then $u_n$ rises only exponentially, rather than factorially, in $n$. When $a-b\in-\N$ we obtain the same deduction from the Kummer symmetry.



\section{Special cases}

In this section we deal with a variety of cases that arise in the literature, and we also deal with all the examples that we presented in the Introduction. The list is not exhaustive, of course, and is simply designed to give the reader a feel for how the method works.

\subsection{Factorial class $S(1,\alpha-1,0)$}

The solution is clearly 
\[
u_n=\frac{\Gamma(n+\alpha)}{\Gamma(1+\alpha)}
\]
which is understood to terminate whenever $-\alpha\in\N$ (and $\myfunc$ is a polynomial of degree $-\alpha$ then).
The recurrence is linear and the Riccati equation is simply
\[
\myfunc'(z) = \left(1+\frac{\alpha}{z}\right)\myfunc(z) - \frac{1}{\alpha}
\]
(where we temporarily restrict attention to $\alpha>0$). The Riccati equation is solved by
\[
\myfunc(z) = \frac{1}{\alpha}ze^zE_{\alpha}(z)
\]
with $E_\alpha$ denoting the exponential integral function \cite[\S5]{Abramowitz64}.
Using the analytic continuation of $E_\alpha$ to the left half-plane, we can write down
\[
\mu(x) = \frac{1}{2\pi\I} [\myfunc]_-(x) = \frac{1}{\Gamma(\alpha+1)} x^{\alpha}e^{-x}
\]
and thereby ascertain that the $(u_n)$ are moments of the Gamma distribution with shape parameter $\alpha$ (note now that the result becomes valid for $\alpha>-1$).

When $\alpha<-1$ we obtain a sequence that is not simply represented, but almost simply represented: $\myfunc(z)$ is of order $z^{\alpha}$ at the origin and $z\myfunc(z)$ is not bounded.
In the special case where $\alpha=-m$, $m\in\N$, the sequence `terminates' at $u_m$ and subsequent terms are all zero, and then $\myfunc(z)$ is of the form $z^{-m}$ multiplied by a polynomial.

This seems like a very heavy-handed way to prove an obvious result (and even more so if the CHGF is invoked) but, as we are about to see, we have not yet finished with the exponential integral function.

\subsection{Recipro-factorial class $S(1,-1-\alpha,\alpha)$}

By reciprocating the previous sequence we can obtain more results. Indeed,
\[
\myfunc(z) = \frac{- e^{-z}}{\alpha z E_{\alpha}(z)}
\]
or if we take the CHGF route, we have $\mathcal{R}(-\alpha,-1,-\alpha,0)$ and $a=b=\alpha$, $k=1$, which gives
\[
w(z) = e^{-z} U(\alpha,\alpha,z) = z^{1-\alpha} E_{\alpha}(z).
\]
Then
\begin{eqnarray*}
U_R(x) &=& e^{-x} \!\left( \Gamma(1-\alpha) + (\cos \pi \alpha)\, x^{1-\alpha} \sum_{r=0}^\infty \frac{1}{r+1-\alpha} \frac{x^r}{r!}  \right) \\
U_I(x) &=&  -(\sin\pi\alpha)\, x^{1-\alpha} e^{-x} \sum_{r=0}^\infty \frac{1}{r+1-\alpha} \frac{x^r}{r!}.
\end{eqnarray*}
By Lemma \ref{lem:zc1}, there are no zeros of $E_\alpha$ in the cut plane for $\alpha>-1$. (The case $\alpha=0$ is unimportant; it is apparently singular because $u_0=\frac{1}{0}$, but this can be removed by subtracting the constant $-\frac{1}{\alpha}$ from $\myfunc(z)$ before allowing $\alpha\to0$, and then one obtains $\myfunc(z)=-e^zE_1(z)$, with $\mu(x)=e^{-x}$.)

One subcase that simplifies a little is when $\alpha=\half+m$, $m\in\N_0$. Then, the $\cos\pi\alpha$  term drops out and we are left with
\[
U_R(x) = e^{-x} \Gamma(\shalf-m), \quad 
U_I(x) = (-)^{m+1} x^{\half-m} e^{-x} \sum_{r=0}^\infty \frac{1}{r+\half-m} \frac{x^r}{r!}.
\]
For example, $m=0$ gives $S(1,-\frac{3}{2},\half)$ which on scaling by 2 gives the known result for $S(2,-3,1)$ stated in the Introduction. (There is no singularity at the origin for $m=0$.  For $m\ge 1$ one has $r^\circ=\frac{1-2m}{1+2m}$.)

When $\alpha\in \N$ the expressions need reworking  as the $\Gamma(1-\alpha)$ term and the summation are both singular. The resulting expressions are
\[
U_R(x) = \frac{(-)^{m}}{(m-1)!} e^{-x} \Ei_{m}(x) , \quad
U_I(x) = \frac{(-)^{m}}{(m-1)!} \pi e^{-x} , 
\qquad \alpha=m\in\N,
\]
with 
\[
\Ei_{m}(x) = \ln x - \digamma(m) + (m-1)! \!\!\!\!\!\!\!\!\sum_{r=0;\,r\ne m-1}^\infty \frac{x^{r+1-m}}{(r+1-m)r!}  .
\]
Now $\myfunc$ is singular at the origin and, as $E_{m}(0)=\frac{1}{m-1}$, $m\ge1$, a simple pole of residue $\frac{1-m}{m}$ is generated. 
In the special case  $m=1$ there is no extra contribution and we just have
\[
U_R(x) = e^{-x} \!\left( -\gamma -\ln x - \sum_{r=1}^\infty \frac{x^r}{r.r!}\right) \equiv -e^{-x}\Ei(x), \qquad U_I(x) = -\pi e^{-x},
\]
(note $\Ei_1\equiv \Ei$), which gives
\[
S(1,-2,1) : \quad \mu(x) = \frac{x^{-1}e^x}{\Ei(x)^2+\pi^2},
\]
as stated earlier. More generally for $\alpha=m\in\N$ we have
\begin{equation}
S(1,-1-m,m) : \quad \mu(x) = \frac{x^{-m}e^x (m-1)!/m}{\Ei_{m}(x)^2+\pi^2} + \textstyle \frac{m-1}{m}\delta(x).
\end{equation}
These sequences are classified in the OEIS (though the closed-form solutions are not; they are purely expressed through the asymptotic expansion of the Riccati equation) but except for $m=1$ they do not seem to have any special combinatorial significance.

Finally for $\alpha\in-\N_0$ we have another important case, as $w(z)$ is then a polynomial $\times \;e^{-z}$. We deal with this later when considering a larger class that  we are calling the \emph{Laguerre} class.

\subsection{Airy cases $S(6,-8,1)$ and $S(6,-4,-1)$}

For $S(6,-8,1)$ we have the Riccati equation as $\mathcal{R}(-\frac{1}{6},0,-\frac{1}{3},\frac{1}{24})$, so the relevant CHGF solution is parametrised by $a=\frac{1}{6}$, $b=\frac{1}{3}$, $k=\frac{1}{6}$. The solution is thereby given via
\[
w(z) = e^{-z/12}U(\textstyle\frac{1}{6},\frac{1}{3},\frac{1}{6}z) = 2^{2/3} 3^{1/6} \pi^{1/2} \Ai \big(\frac{1}{4}z^{2/3}\big)
\]
or equivalently
\[
\myfunc(z) = \frac{z^{-1/3}\Ai'(\frac{1}{4}z^{2/3})}{\Ai(\frac{1}{4}z^{2/3})}
\]
Now for $x\in\R^+$ we have
\begin{eqnarray*}
\Uplus(\textstyle\frac{1}{6},\frac{1}{3},-\frac{1}{6}x) &=& e^{-x/12} 2^{2/3}3^{1/6}\pi^{1/2} \Ai\big(\textstyle\frac{1}{4}\omega x^{2/3}\big)\\
\Uminus(\textstyle\frac{1}{6},\frac{1}{3},-\frac{1}{6}x) &=& e^{-x/12} 2^{2/3}3^{1/6}\pi^{1/2} \Ai\big(\textstyle\frac{1}{4}\baromega x^{2/3}\big),
\end{eqnarray*}
with $\omega=e^{2\pi\I/3}$ and $\baromega=e^{-2\pi\I/3}$.
Using the identities
\[
\Ai (e^{\pm2\pi\I/3} x) = \shalf e^{\pm \pi\I/3} \big(\Ai(x) \mp \I \Bi(x)\big)
\]
we find
\[
\Uplus(\textstyle\frac{1}{6},\frac{1}{3},-\frac{1}{6}x) \Uminus(\textstyle\frac{1}{6},\frac{1}{3},-\frac{1}{6}x) = 2^{-2/3}3^{1/3}\pi e^{-x/6} \big[ \Ai(\frac{1}{4}x^{2/3})^2 + \Bi(\frac{1}{4}x^{2/3})^2\big]
\]
whence by (\ref{eq:result1}) we obtain
\begin{equation}
S(6,-8,1): \quad \mu(x)=\frac{\pi^{-2}x^{-1/3}}{\Ai(\frac{1}{4}x^{2/3})^2 + \Bi(\frac{1}{4}x^{2/3})^2},
\end{equation}
as claimed earlier. The sequence is regular because $\Ai(\frac{1}{4}z^{2/3})$ has no zeros for $|\arg z|<\frac{3}{2}\pi$ (a known result \cite{Abramowitz64}).

The sequence $S(6,-4,-1)$ is its reciprocal, and the corresponding Riccati equation is $\mathcal{R}(\frac{1}{6},0,\frac{1}{3},-\frac{1}{24})$ whose solution is
\[
\myfunc(z) = -\frac{z^{1/3}\Ai(\frac{1}{4}z^{2/3})}{4\Ai'(\frac{1}{4}z^{2/3})}
\]
which is generated from
\[
\textstyle
w(z) = e^{-z/12}U(-\frac{1}{6},-\frac{1}{3},\frac{1}{6}z) 
= -\pi^{1/2} 2^{4/3} 3^{-1/6} \Ai'(\frac{1}{4}z^{2/3}).
\]
Hence
by (\ref{eq:result1}) we obtain, as $\Ai'$ also has no zeros in the cut plane,
\begin{equation}
S(6,-4,-1): \quad \mu(x)=\frac{\frac{1}{4}\pi^{-2}x^{1/3}}{\Ai'(\frac{1}{4}x^{2/3})^2 + \Bi'(\frac{1}{4}x^{2/3})^2},
\end{equation}
which begins
\[
S(6,-4,-1) = (1,7,84,1463,33936,990542, \ldots).
\]
The application of this is that it determines the moments of the area of the absolute value of a Brownian bridge \cite{Janson07,Kearney09}.

As a brief demonstration of the asymptotic analysis that we derived earlier, we have
the following table. For both sequences we have $u_n \sim 6^n(n-1)!/2\pi$.

\vspace{2mm}

\begin{center}\begin{tabular}{rrrr}
\hline
$n$ & $\ln S(6,-8,1)_n$ & $\ln S(6,-4,-1)_n$ & Asymptotic \\
\hline
1 & 0.000 & 0.000	& $-$0.046 \\
2 & 1.609 & 1.946 &	1.746\\
3 & 4.094 & 4.431 &	4.231\\
4 & 7.008 & 7.288 &	7.121\\
5 & 10.208	& 10.432 & 10.299\\
6 & 13.627	& 13.806  & 13.700\\
7 & 17.224	& 17.369 & 17.284\\
8 & 20.972	& 21.092 & 21.021\\
9 & 24.850	& 24.952 & 24.893\\
10 & 28.845	& 28.933 & 28.882\\
\hline
\end{tabular}\end{center}

\vspace{2mm}\noindent
In this case the relative error in estimating $u_n$ is ca.~5\% for $n\ge 10$.

\subsection{Bessel class $S(1,2\nu-2,\half-\nu)$, $-\frac{3}{2}\le\nu\le\frac{3}{2}$}

This is a generalisation of the Airy example that we have just treated; they have in common the property that $\beta_2=0$, or equivalently $\alpha_3=-\half(\alpha_1+\alpha_2)$.
The Riccati equation is $\mathcal{R}(\beta,0,2\beta,-1/4\beta)$ with $\beta=\nu-\half$, and the CHGF parameters are $a=\half-\nu$, $b=2a$.
The Airy cases taken above are given by $\nu=\frac{1}{3},\frac{2}{3}$. The sequence is irregular for $\nu$ outside the range $[-\frac{3}{2},\frac{3}{2}]$, and by 
Lemma \ref{lem:zc1} there are no zeros in the cut plane for $0\le \nu \le\frac{3}{2}$.
 The class is closed under the Kummer transformation, which simply effects the transformation $\nu\mapsto-\nu$, so we need only analyse $0\le \nu \le \frac{3}{2}$, which is regular. The Bessel class is also self-reciprocal: the reciprocal sequence is parametrised by $\nu^\dagger = 1-\nu$. Note that the reciprocal may not preserve regularity, e.g.\ $\nu=-1$ gives a quasiregular sequence but $\nu=1-(-1)=2$ gives an irregular one.

In the Bessel case, the solution is
\[
w(z) = e^{-z/2} U(\shalf-\nu,1-2\nu,z) = \pi^{-1/2} z^\nu K_\nu(\shalf z)
\]
(note that $K_\nu\equiv K_{-\nu}$), and
\[
\myfunc(z) = \frac{1}{2\nu-1} \frac{K_{\nu-1}(\half z)}{K_\nu(\half  z)}
\]
The self-reciprocity can be seen from this solution, as $(z^\nu K_\nu)'=-z^\nu K_{\nu-1}$ and $(z^{-\nu} K_\nu)'=-z^{-\nu}K_{\nu+1}$.

The cases $\nu=0,1$ may therefore be handled together and can be treated easily enough using the expansions of $K_0$, $K_1$ and their properties. We have
\[
\myfunc(z) = -\frac{K_1(\half z)}{K_0(\half z)} \mbox { or } \frac{K_0(\half z)}{K_1(\half z)} 
\]
respectively (these can also be obtained from the relations $K_0'=-K_1$, $(zK_1)'=-zK_0$).
The analytic continuation of $K_0(z)$ round to the left half-plane is
\[
K_0^\oplus(-x) = K_0(x) - \pi \I  I_0(x), \quad K_0^\ominus(-x) = K_0(x) + \pi \I I_0(x),
\]
so we obtain
\begin{eqnarray}
S(1,-2,\shalf)&:& \quad \mu(x)=\frac{2/x}{K_0(\half x)^2 + \pi^2 I_0(\half x)^2} \\
S(1,0,-\shalf)&:& \quad \mu(x)=\frac{2/x}{K_1(\half x)^2 + \pi^2 I_1(\half x)^2} \nonumber.
\end{eqnarray}
On scaling by 2 one obtains the results mentioned in the Introduction.

When $\nu$ is not an integer a direct appeal to the CHGF has to be made, and we evaluate $U_R$ and $U_I$ using the result
\[
M(\shalf+\nu,1+2\nu,x) = 4^{\nu}\Gamma(1+\nu) x^{-\nu}e^{x/2} I_\nu(\shalf x) .
\]
The duplication and reflection formulae for the Gamma function allow all the $\Gamma$ terms to be cleared up, and the result is
\begin{eqnarray*}
U_R(x) &=& \pi^{-1/2} x^\nu e^{-x/2} \left(  K_\nu(\shalf x) + \pi (\sin\pi\nu)\,  I_{\nu}(\shalf x)\right) \\
U_I(x) &=& -\pi^{1/2}(\cos\pi\nu ) \, x^\nu e^{-x/2} I_\nu(\shalf x),
\end{eqnarray*}
where we have also used the expression for $K_\nu$ in terms of $I_{\pm\nu}$,
\[
K_\nu(z) = \frac{\pi}{2\sin\pi\nu} \big( I_{-\nu}(x)-I_\nu(x)\big).
\]
The above expression for $U_R,U_I$ is also valid for integer $\nu$. 

By (\ref{eq:result1}) we obtain
\begin{equation}
S(1,2\nu-2,\shalf-\nu): \quad \mu(x)=\frac{(\half-\nu)^{-1} (\cos\pi\nu ) /x}
{K_\nu(\half x)^2 + \pi^2I_\nu(\half x)^2 + 2 \pi(\sin\pi\nu)\, K_\nu(\half x) I_\nu(\half x) }.
\label{eq:bessmu}
\end{equation}
As expected, the Airy case can be obtained from this by setting $\nu=\frac{1}{3},\frac{2}{3}$.
Now the cases $\nu=\frac{1}{2},\frac{3}{2}$ cause minor difficulties in (\ref{eq:bessmu}). However, $\mu(x)$ just boils down to $e^{-x}$ in the case $\nu=\half$, and so is not of particular interest given that we have already dealt with the sequence $S(1,-1,0)$. When $\nu=\frac{3}{2}$ we have $S(1,1,-1)$ which is simply the powers of 2. This is because $\myfunc(z)$ has a pole on the negative real axis (at $z=-1$), and indeed the denominator of (\ref{eq:bessmu}) completes as a square that vanishes at $x=2$, so $\mu(x)$ is to be interpreted as $\delta(x-2)$.

\subsection{Miscellaneous}

The sequence $S(2,-2,1)$, which is self-reciprocal, has $a=\half$, $b=0$, $k=\half$. From Lemmas \ref{lem:origin} and \ref{lem:zc1}, the problem is regular.
The expression for $U_R(x)$ is singular at $b=0$ by virtue of the Gamma functions, but the singularity is removable and application of L'H\^opital's rule yields
\[
U_R(x) = e^{-x} \left[ \frac{2}{\sqrt{\pi}} - \frac{\ln x}{\pi} \sum_{r=0}^\infty \frac{\Gamma(\half+r)}{\Gamma(2+r)} \frac{x^{r+1}}{r!} + \frac{1}{\pi} \sum_{r=0}^\infty \frac{\Gamma(\half+r)}{\Gamma(2+r)} d_r \frac{x^{r+1}}{r!}\right]
\]
with
\[
d_r = \digamma(2+r)-\digamma(\shalf+r)+\digamma(1+r);
\]
and
\[
U_I(x) = -xe^{-x} \sum_{r=0}^\infty  \frac{\Gamma(\half+r)}{\Gamma(2+r)} \frac{x^r}{r!}.
\]
After some minor manipulations one ends up with the result given in the Introduction.

The related sequence (\#A005413)
\[
(1, 7, 72, 891, 12672, 202770, \ldots),
\]
which is obtained by $v_n=(n-1)(u_n+2nu_{n-1})$, $n\ge 2$, is thereby solved too, and in QED is the number of graphs with `proper vertices' \cite{Cvitanovic78}.

\subsection{Laguerre class $S(-1,2+\alpha,m)$, $S(-1,2-\alpha,\alpha+m)$, $m\in\N$}

The second of these cases is the Kummer transformation of the first, so we only have to deal with the first. In this case $k=-1$, $a=-m$, $b=1+\alpha$ and then the CHGF is simply a polynomial. (We have chosen to reverse the signs and have $k$ negative so as to keep the branch cut along the negative real axis.) The solution is given by
\[
w(z) = e^{-mz/(1+\alpha)} U(-m,1+\alpha,-z)
\]
which is expressible in terms of the Laguerre polynomial defined by
\[
L^{(\alpha)}_m (z) = \frac{(-)^m}{m!}U(-m,1+\alpha,z) = \frac{1}{2\pi\I} \oint_{(0+)} e^{-zt} t^{-m-1}(1+t)^{m+\alpha}\, dt.
\]
One therefore has
\[
\myfunc(z) =   \frac{L^{(\alpha)}_m {}' (-z)}{m\, L^{(\alpha)}_m (-z)} + \frac{1}{1+\alpha}
\]
which admits a partial fraction expansion
\[
\myfunc(z) =   \frac{1}{1+\alpha} + \frac{1}{m} \sum_{j=1}^m \frac{1}{\zeta_j-z}.
\]
Intriguingly, the sequence therefore obeys a \emph{linear} recurrence relation of order $m$, with constant coefficients. At this point one might jump to the conclusion that as the zeros of the Laguerre polynomials are well-studied, everything must be plain sailing. That is not so.

When $\alpha>-1$ the zeros of $L^{(\alpha)}_m(x)$ are known to be real and positive, as they are orthogonal polynomials of the Gamma distribution (density $\propto x^\alpha e^{-x}$). The theory of these is extensive. For example, the zeros of $L^{(\alpha)}_m(x)$ interlace those of $L^{(\alpha)}_{m+1}(x)$; and asymptotic results are known for large $m$.
The sequence is therefore regular in this case, and  one has, on writing $\xi_i$ for the zeros of $L^{(\alpha)}_m(x)$,
\begin{equation}
S(-1,2+\alpha,m)_n = \frac{1}{m} \sum_{i=1}^m \xi_i^{n-1} \qquad (\alpha>0).
\label{eq:laguerre2}
\end{equation}
So the sequence has a  Mellin representation as a `comb' of delta-functions
\[
\mu(x) = \frac{1}{m} \sum_{i=1}^m \delta(x-\xi_i).
\]
This is a discrete probability measure bearing some resemblance to the Gaussian quadrature formula of the associated Gamma distribution: the locations of the delta-functions are the same, but the strengths of the delta-functions are all equal to $\frac{1}{m}$, unlike the ordinates of the Gaussian quadrature formula.
As an example, $\alpha=0$, $m=2$ gives
\[
S(-1,2,2)_n = \textstyle \frac{1}{2}  \big( (2+\sqrt{2})^{n-1} + (2-\sqrt{2})^{n-1}\big) ,
\]
so the delta-functions are situated at $x=2\pm\sqrt{2}$. It is readily verified that $u_n=4u_{n-1}-2u_{n-2}$.

But when $\alpha\le -1$ the connection with Gaussian quadratures disappears, and now the zeros need not be real. In general, these sequences are not regular. The relation (\ref{eq:laguerre2}) is still valid but there is in general no representation as a measure.
As examples, with $\alpha=-1$ and $m=3$ we get
\[
S(-1,1,3)_n = \textstyle \frac{1}{3}  \big( (3+\sqrt{3})^{n-1} + (3-\sqrt{3})^{n-1} + 0^{n-1} \big) ,
\]
which is quasiregular. On the other hand $\alpha=-4$, $m=3$ gives
\[
S(-1,-2,3) = ( 1,-1,-1,3,3,-21,27,27,-117,27,459,-837,\ldots )
\]
which is irregular (and in appearance too) and solved by (\ref{eq:laguerre2}) with roots given as
\[
\xi_1 = \textstyle\sqrt[3]{\!\!\sqrt{2}-1} - \sqrt[3]{\!\!\sqrt{2}+1} - 1 \\
\]
\[
\xi_{2,3} = \frac{\textstyle\sqrt[3]{\!\!\sqrt{2}+1} - \sqrt[3]{\!\!\sqrt{2}-1}}{2} -1 \pm \frac{\I\sqrt{3}}{2} \big( \textstyle\sqrt[3]{\!\!\sqrt{2}-1} + \sqrt[3]{\!\!\sqrt{2}+1}\big) ,
\]
these being the roots of $x^3+3x^2+6x+6=0$. The associated linear recurrence is $u_n=-3u_{n-1}-6u_{n-2}-6u_{n-3}$. Note incidentally that up to scaling by $-1$ this is $S(1,2,-3)$ which is in the recipro-factorial class.

One thing that does carry over from the theory of Gaussian quadratures irrespective of $\alpha$ is the connection with Jacobi, or tridiagonal, matrices \cite{Deift00}. By exhibiting an $m$-dimensional tridiagonal matrix $T$ whose characteristic polynomial is the monic Laguerre polynomial, we can simply write (\ref{eq:laguerre2})
as
\begin{equation}
u_n = \frac{1}{m}\, \tr [T^{n-1}],
\end{equation}
where $\tr$ denotes the trace and $T^0$ is the identity matrix. An explicit expression for $T$ is
\begin{equation}
T = \begin{bmatrix} 1+\alpha & 1 &  \\ 1(1+\alpha) & 3+\alpha & 1 \\  & 2(2+\alpha) & 5+\alpha & 1 \\
&& \ddots & \ddots &   \\
&&& m'(m'+\alpha) & 2m'+1+\alpha 
\end{bmatrix},
\end{equation}
where $m'=m-1$.

\section{Algebraic class $S(0,\alpha,1)$}
\label{sec:Alg}

We round off this paper by treating what is on the face of it a rather trivial special case: when $\alpha_1=0$, the Riccati equation reduces to a quadratic (algrebraic) equation. It turns out that this case is distinct, one reason being that different types of singularities are encountered from the general case that we solved using the CHGF. 

This case is definitely of practical interest as it crops up frequently in combinatorics, most obviously through the  well-known Catalan numbers,
\[
C_n = S(0,0,1)_{n+1}; \qquad C_{n+1} = \sum_{j=0}^n C_j C_{n-j}.
\]

We now proceed to solve this case for all $\alpha$. The Riccati equation is solved as 
\[
\myfunc(z) = \frac{(1+2/\alpha)z+\alpha - \sqrt{z^2+2(2+\alpha)z+\alpha^2}}{2z}
\]
which needs careful interpretation, as follows.

\begin{itemize}
\item $\alpha>0$. 
The plane is cut between $z=-2-\alpha\pm2\sqrt{1+\alpha}$ and the square root in the above equation is taken as real and positive when $z$ is. At $z=0$ the numerator of the above expression for $\myfunc$ is $\alpha-\sqrt{\alpha^2}=0$, so $\myfunc$ is regular there.  By considering the jump in $\myfunc$ across the branch cut, we identify the series as regular and with a measure that is \emph{bounded}:
\begin{equation}
\mu(x) = \left\{ \begin{array}{ll}  \displaystyle  \frac{\sqrt{-x^2+2(2+\alpha)x-\alpha^2}}{2\pi x}, & |x-(2+\alpha)| \le 2\sqrt{1+\alpha} \\ 0, & \mbox{otherwise.} \end{array} \right.
\label{eq:randmat1}
\end{equation}
\item $\alpha=0$ (Catalan). We can let $\alpha\to0$ in the above case without vitiating the result, despite the singular nature of $\myfunc$ (which is caused by $u_0=\frac{1}{\alpha_1+\alpha_2}=\infty$ and hence is of no consequence). This gives
\begin{equation}
\mu(x) = \left\{ \begin{array}{ll}  \displaystyle  \frac{\sqrt{4-x}}{2\pi \sqrt{x}}, & |x-2| < 2 \\ 0, & \mbox{otherwise.} \end{array} \right.
\label{eq:randmat2}
\end{equation}
\item $-1<\alpha<0$. The plane is cut as before with the same definition of the square root. But now $\alpha-\sqrt{\alpha^2}=2\alpha$, and so $\myfunc$ has a simple pole at the origin of residue $\alpha$. Thus a delta-function is created:
\begin{equation}
\mu(x) = \left\{ \begin{array}{ll}  \displaystyle  \frac{\sqrt{-x^2+2(2+\alpha)x-\alpha^2}}{2\pi x}  , & |x-(2+\alpha)| < 2\sqrt{1+\alpha} \\ (-\alpha) \delta(x), & \mbox{otherwise.} \end{array} \right.
\label{eq:randmat3}
\end{equation}
\item $\alpha=-1$. Here $\myfunc(z)=-1/z$ and so $\mu(x)=\delta(x)$.  
\item $-2<\alpha<-1$. The branch points are complex with negative real part; again the square root is taken to be real and positive whenever $z$ is.
Write $\tilde{\alpha}=\sqrt{|1+\alpha|}$ and parametrise the complex variable on the branch cut by $\zeta_y \equiv -2-\alpha+2\I \tilde{\alpha} y$. We arrive at
\[
u_n =  -\frac{2\tilde{\alpha}^2}{\pi}  \int_{-1}^1 (-\zeta_y)^{n-2} \sqrt{1-y^2}  \, dy , \qquad n>1.
\]
The contribution from the pole at the origin takes care of the $u_1$ term.
\item $\alpha=-2$. Here $\myfunc(z)=-(1+\sqrt{z^2/4+1})/z$. The even-numbered terms are the Catalan numbers, with alternating signs, and the odd terms vanish after the first:
\[
S(0,-2,1)_{2n} = (-)^n C_{n-1} = -\frac{2}{\pi} \int_{-1}^1 (-4y^2)^{n-1} \sqrt{1-y^2} \, dy
\]
This can also be seen from a direct appeal to the Binomial expansion of $\myfunc(z)$.
\item $\alpha<-2$. This case is similar to $-2<\alpha<-1$, but the branch points are in the right half-plane and the origin lies to the left of them, so $\myfunc$ is now \emph{regular} at the origin. We therefore obtain, as before,
\[
u_n = -\frac{2\tilde{\alpha}^2}{\pi}  \int_{-1}^1 (-\zeta_y)^{n-2} \sqrt{1-y^2}  \, dy, \qquad n\ge 1.
\]

\end{itemize}

A convenient summary is:
\begin{thm}
The sequence $S(0,\alpha,1)$ is regular if $\alpha\ge0$, and quasiregular if $-1\le\alpha<0$. If $\alpha<-1$, then it is not simply represented. $\Box$
\end{thm}

It appears that there is an interesting connection between some of these cases (\ref{eq:randmat1}--\ref{eq:randmat3}) and the theory of random matrices, and we pursue this discussion in the last section of the paper.

\subsubsection*{Trace formula}

In all cases for $\alpha \ge-1$ the sequence is regular and admits a trace formula in the usual way.
The other situations do not quite fit the pattern we have previously encountered.
For $\alpha<-1$, the integral
\[
u_n =  -\frac{2\tilde{\alpha}^2}{\pi}  \int_{-1}^1 (-\zeta_y)^{n-2} \sqrt{1-y^2}  \, dy 
\]
bears a resemblance to a continuous version of the expression $\frac{1}{m}\sum_{j=1}^m (-\zeta_j)^{n-1}$ that we encountered in the irregular Laguerre case, and it is tempting to think of it as an integral over some `pole density' $\widehat{\mu}(y)$. 
However, the measure would have to be
\[
\widehat{\mu}(y) := -\left(\frac{2\tilde{\alpha}^2}{\pi}\right) \frac{\sqrt{1-y^2}}{2+\alpha-2\I \tilde{\alpha} y}, \qquad |y|<1,
\]
which is not real, let alone real and positive. So there is no obvious interpretation as a limiting distribution of poles arising from the Laguerre case. An open question is, therefore, how to interpret the complex density, possibly by trying to split it into real and imaginary parts and examining each separately.
\notthis{
That said, some manipulation gives
\[
u_n = - \frac{2\tilde{\alpha}^2}{\pi} \int_{-1}^1 (-\zeta_y)^{n-1} \frac{(2+\alpha)\sqrt{1-y^2}}{(2+\alpha)^2+\tilde{\alpha}^2 y^2} \, dy
- \frac{2\tilde{\alpha}^2}{\pi} \int_{-1}^1 (-\zeta_y)^{n-1} \frac{2\I\tilde{\alpha} y \sqrt{1-y^2}}{(2+\alpha)^2+\tilde{\alpha}^2 y^2} \, dy
\]} 



\section{Conclusions and Discussion}

We have considered the sequence (\ref{eq:u_n}) and shown how it may be solved in closed form as the expression (\ref{eq:gensoln}). In the case we have described as \emph{regular} the solution reduces to (\ref{eq:Mellin}).

We list some questions and applications that naturally arise from this work.

\subsubsection*{Trace formula}

It seems that the expression for $u_n$ is the trace of the $(n-1)$th power of a linear operator. In the Laguerre case the spectrum of this operator is discrete and it admits a matrix representation. Can this representation as the trace of an iterated operator be extended to all the sequences? It seems reasonable to suppose that the operator has something to do with the physical system that gives rise to the sequence.

\subsubsection*{Random matrices}

It appears that there is an interesting connection between some of these cases (\ref{eq:randmat1}--\ref{eq:randmat3}) and the theory of random matrices.
 Let $X$ be a $(q\times p)$ matrix whose entries are independent and Normally distributed   with mean zero and variance 1. Consider now the empirical covariance matrix
\[
\Sigma = \frac{1}{\max(p,q)} XX',
\]
in which $'$ denotes transpose, and  let $p,q\to\infty$ with $p/q$ held fixed. The limiting distribution of the eigenvalues of this $(q\times q)$ matrix is, provided $p>q$, 
\[
\frac{1}{2\pi(q/p)} \frac{\sqrt{(\lambda_+-x)(x-\lambda_-)}}{x} , \qquad \lambda_\pm = (1\pm\sqrt{q/p})^2,
\]
which up to a scaling is identical to $\mu(x)$ in eq.~(\ref{eq:randmat1}). This result was derived rigorously by V.~Mar{\v c}enko and L.~Pastur \cite{Marcenko67} and studied more recently by A.~Sengupta and P.~Mitra \cite{Sengupta99}. 

When $-1<\alpha<0$, the connection still works. This is the case when $p<q$, so $\Sigma$ is rank-deficient and must have at least $q-p$ zero eigenvalues: hence the occurrence of a point probability mass at zero in (\ref{eq:randmat3}). Of course the case $p=q$ is covered by (\ref{eq:randmat2}), which on changing variables $x \rightarrow x^2$ states that the square roots of the eigenvalues are distributed according to the `semicircle law' $\sqrt{4-x^2}/\pi$ $(0<x<2)$, a famous result of E.~Wigner \cite{Wigner55,Wigner58}.

Both of these can be understood as limiting instances of the Laguerre class $S(-1,\alpha m,m)$ as $m\to\infty$. The Laguerre polynomials also arise in the theory of random matrices, and clearly this is an avenue worth investigating.

\subsubsection*{Numerical analysis: Finding zeros of $U(a,b,z)$}

When the sequence is irregular, $\myfunc$ has poles in the cut plane. So far the emphasis has been on finding a closed form solution to $u_n$, which would necessitate counting them and finding their locations. However, if we use the methods in reverse, and compare the integral $\int_0^\infty x^{n-1}\mu(x)\,dx$ with $u_n$ for a few values of $n$, we can find them. An example is adequate to demonstrate the general case, and we take the Bessel case with $\nu=2$.
We have $a=-\frac{3}{2}$, $b=-3$, and there is no contribution from the origin. By Lemma \ref{lem:zc1} there are two zeros to find. We therefore have
\[
S(1,2,{-\textstyle\frac{3}{2}})_n =   -{\textstyle \frac{2}{3}}
\int_0^\infty  \frac{x^{n-2}}{K_2(\half x)^2 + \pi^2I_2(\half x)^2} \, dx
+ {\textstyle \frac{2}{3}} \sum_{j=1}^2 (-\zeta_j)^{n-1},
\]
where $K_2(\frac{1}{2}\zeta_j)=0$.
Substituting $n=1$ simply confirms that we were right about the number of zeros, as the integral (without the prefactor) is equal to $\frac{1}{2}$.
\notthis{
For large $x$ the integrand is asymptotically $e^{-x}$ multiplied by a constant ***which is what?***, so this directs us to calculate the integral by Gauss-Laguerre quadrature, i.e.\
\[
\int_0^\infty  \left[\frac{x^{n-2}e^x}{K_2(\half x)^2 + \pi^2I_2(\half x)^2}\right] e^{-x} \, dx.
\]
The term in square brackets, implicitly being approximated by a polynomial; as it is well-behaved, we expect this quadrature method to work efficiently.
} 
For higher $n$ we evaluate the integral  numerically, which presents no difficulties as the integrand is smooth and unimodal. We arrive at the following table:

\vspace{3mm}

\begin{center}\begin{tabular}{rrrr}
\hline
$n$ & $S(1,2,{-\textstyle\frac{3}{2}})_n$ & Integral & Difference, $v_n$ \\
\hline
1 & 1 & $-0.333333$ & 1.333333 \\
2 & 2.5 & $-0.916997$ & 3.416997 \\
3 & 5 & $-2.77313$ & 7.77313 \\
4 & 5.625 & $-9.25324$ & 14.87824 \\
5 & $-15$ & $-34.4717$ & 19.4717 \\
6 & $-154.6875$ & $-145.797$ & $-$8.89074 \\
\hline
\end{tabular}\end{center}

\vspace{3mm}
\noindent
The zeros obviously occur as a complex-conjugate pair, $re^{\pm \I\theta}$ say. Therefore the difference is equal to
\[
v_n = \textstyle \frac{4}{3} r^{n-1} \cos((n-1)(\pi-\theta)).
\]
We now employ Prony's method (see e.g.~\cite{Marple87}) to extract $r,\theta$. This consists in observing that $(v_n)$ obeys a linear recurrence, here of order 2:
\[
v_n + c_1v_{n-1} + c_2v_{n-2}=0,
\]
with the coefficients identifiable easily enough from $(v_n)$ by linear regression, after which the $(-\zeta_j)$ are the roots of the quadratic equation $z^2+c_1z+c_2=0$ (equivalently, $c_1=2r\cos\theta$ and $c_2=r^2$). Applying this to the above,  
we find
$c_1=-5.12549$, $c_2=7.305500$, 
from which we deduce that the zeros of $K_2(\half z)$ in the cut plane are located at
\[
\zeta \approx 2.70287 e^{\pm2.81818\I}.
\]
A few points can be made about the above analysis. First, we have identified the \emph{complex} zeros of $K_2(\half z)$ \emph{despite evaluating the function at real values only}, which has a certain appeal to it and is particularly useful if one does not have at one's disposal an algorithm for its evaluation at complex arguments. The method is applicable to higher orders of modified Bessel function, which have more zeros. In practice, we suggest that the best way to find zeros is probably to use the technique to obtain reasonable accuracy and then use Newton's method to polish them one by one. On the zeros of $K_2(\half z)$, using the above estimate as starting-point, we find the positions of the roots (using the computer package MATLAB) to be
in agreement with the above, to 6 s.f.
This method has an important advantage once the number of zeros becomes reasonably large. Without such assistance, one would have to search for the zeros of a special function without a clear idea of where to find them (Newton's method not being globally convergent in general). The associated polynomial in Prony's method is of degree $\# U$, so for $\# U\le 4$, it can be factorised by radicals; for $\# U >4$, search algorithms have to be used to factorise it, but this is a reasonably standard problem \cite{NRC}, and certainly much easier than finding the zeros of a special function.

\subsubsection*{Number-theoretic properties}
So far, the analysis of the recurrence has been done using complex variable theory. However, there are certain questions that are likely to be resolved more easily by other methods. For example, for how many $n$ does $u_n=5$ (a question that is only interesting for the irregular series in which terms vary in sign, as in $S(-1,-2,3)$ earlier); or for now many $n$ does $u_n=3$ (mod 7)?
Questions of this form are much easier to tackle using $p$-adic analysis \cite{Cassels86,Schikhof84} than complex analysis: one such theorem that is arises from local field theory is the theorem of Mahler \& Lech, which states that if a sequence obeys a linear recurrence with constant coefficients, then the set of $n$ for which $u_n=0$ (or some other specified number) is a finite combination of sets each of which is either periodic or finite. Obviously this theorem does not apply to the recurrence discussed herein, except for the sequences of the Laguerre class, which do obey a linear recurrence. However, the hypergeometric functions used in the paper do have $p$-adic analogues, even though the convergence properties are generally very different in the $p$-adic valuation (essentially because $|n!|_p\to 0$ as $n\to\infty$), so the use of $p$-adic analysis may be fruitful.

\subsubsection*{Other recurrences}
The methods used in this paper could possibly be applied to other self-convolutive recurrences of a more general type. It seems that quantum field theory yields a rich source of integer sequences obeying similar self-convolutive recurrences. 
For example, L.~Molinari \cite{Molinari05} gives some sequences related to the enumeration of Feynman diagrams that arise from asymptotic expansion of solutions of first-order nonlinear differential equations that have a cubic nonlinearity,
\[
\myfunc(z) \big(1+z^{-1}\myfunc(z)\big)^2 = 1-\frac{\myfunc(z)}{z}  +2 \frac{d\myfunc(z)}{dz} ,
\]
(we have altered his notation so that the expansion is in ascending powers of $-z^{-1}$), and hence require further manipulation before they can be reduced to a Riccati equation. A recent article by Y.~Pavlyukh and W.~H\"ubner \cite{Pavlyukh07} solves for the generating function but does not give a final reduction to an expression resembling (\ref{eq:gensoln}); hence further work is required in this direction.
There is also the possibility of using more general types of differential equation such as the general hypergeometric equation, and possibly invoking the Painlev\'e transcendents as a more general class of nonlinear differential equations.

That these sequences arise in other branches of mathematics should make this work quite widely applicable. Aside from this there is an elegance about the result that makes the subject quite attractive, and we hope that this, combined with the applications to the problems from which the sequences arise, will give other researchers the enjoyment that we have had from investigating this topic.

\section*{Acknowledgement}
We thank the reviewers for their generous and attentive comments and suggestions for improvement.

\bibliographystyle{plain}
\bibliography{./phd}

\begin{thebibliography}{10}

\bibitem{Abramowitz64}
M.~Abramowitz and I.~A. Stegun.
\newblock {\em Handbook of mathematical functions}.
\newblock Dover, New York, 1964.

\bibitem{Bender78}
C.~M. Bender and S.~A. Orszag.
\newblock {\em Advanced mathematical methods for scientists and engineers}.
\newblock McGraw-Hill, New York, 1978.

\bibitem{Cassels86}
J.~W.~S. Cassels.
\newblock {\em Local Fields}.
\newblock CUP, 1986.

\bibitem{Cvitanovic77}
P.~Cvitanovi\'c.
\newblock Asymptotic estimates and gauge invariance.
\newblock {\em Nuclear Physics}, B127:176--188, 1977.

\bibitem{Cvitanovic78}
P.~Cvitanovi\'c, B.~Lautrup, and R.~B. Pearson.
\newblock Number and weights of {F}eynman diagrams.
\newblock {\em Phys. Rev. D}, 18(6):1939--1949, 1978.

\bibitem{Deift00}
P.~Deift.
\newblock {\em Orthogonal Polynomials and Random Matrices: A Riemann-Hilbert
  Approach}.
\newblock AMS, 2000.

\bibitem{Flajolet09}
P.~Flajolet and R.~Sedgewick.
\newblock {\em Analytic Combinatorics}.
\newblock CUP, 2009.

\bibitem{Janson03}
S.~Janson.
\newblock The {W}iener index of simply generated random trees.
\newblock {\em Random Structures Algorithms}, 22:337--358, 2003.

\bibitem{Janson07}
S.~Janson.
\newblock Brownian excursion area, {W}right's constants in graph enumeration,
  and other {B}rownian areas.
\newblock {\em Prob. Surveys}, 4:80--145, 2007.

\bibitem{Kearney07}
M.~J. Kearney, S.~N. Majumdar, and R.~J. Martin.
\newblock The first-passage area for drifted {B}rownian motion and the moments
  of the {A}iry distribution.
\newblock {\em J. Phys. A: Math. Theor.}, 40:F863--F869, 2007.

\bibitem{Kearney09}
M.~J. Kearney and R.~J. Martin.
\newblock Airy asymptotics: the logarithmic derivative and its reciprocal.
\newblock {\em J. Phys. A: Math. Theor.}, 42:425201, 2009.

\bibitem{Marcenko67}
V.~A. Mar{\v c}enko and L.~A. Pastur.
\newblock Distribution of eigenvalues for some sets of random matrices.
\newblock {\em Math. USSR Sbornik}, 1(4):457--483, 1967.

\bibitem{Marple87}
S.~L. {Marple Jr.}
\newblock {\em Digital spectral analysis with applications}.
\newblock Prentice-Hall, Englewood Cliffs, NJ, 1987.

\bibitem{Molinari05}
L.~G. Molinari.
\newblock Hedin's equations and enumeration of {F}eynman diagrams.
\newblock {\em Phys. Rev.}, B71:113102, 2005.

\bibitem{Pavlyukh07}
Y.~Pavlyukh and W.~H\"ubner.
\newblock Analytic solution of {H}edin's equations in zero dimensions.
\newblock {\em J. Math. Phys.}, 48:052109, 2007.

\bibitem{NRC}
W.~H. Press, B.~P. Flannery, S.~A. Teukolsky, and W.~T. Vetterling.
\newblock {\em Numerical Recipes in C++}.
\newblock CUP, 2002.

\bibitem{Schikhof84}
W.~Schikhof.
\newblock {\em Ultrametric Calculus}.
\newblock CUP, 1984.

\bibitem{Sengupta99}
A.~M. Sengupta and P.~P. Mitra.
\newblock Distributions of singular values for some random matrices.
\newblock {\em Phys. Rev. E}, 60:3389--3392, 1999.
\newblock Also arXiv:cond-mat \#9709283 (1997).

\bibitem{Sloane73}
N.~J.~A. Sloane.
\newblock {\em Handbook of Integer Sequences}.
\newblock Academic Press, 1973.
\newblock Online at {\tt www.research.\-att.com/$\sim$njas/sequences}.

\bibitem{Wigner55}
E.~Wigner.
\newblock Characteristic vectors of bordered matrices with infinite dimensions.
\newblock {\em Ann. Math.}, 62:548--564, 1955.

\bibitem{Wigner58}
E.~Wigner.
\newblock On the distribution of the roots of certain symmetric matrices.
\newblock {\em Ann. Math.}, 67:325--328, 1958.

\end{thebibliography}

\end{document}